\documentclass[11pt]{amsart}
\usepackage{epsfig,amsmath}

\newtheorem{thm}[subsection]{Theorem}
\newtheorem{lem}[subsection]{Lemma}

\newtheorem{prop}[subsection]{Proposition}
\newtheorem{obs}[subsection]{Observation}

\theoremstyle{definition}
\newtheorem{definition}[subsection]{Definition}

\newtheorem{remark}[subsection]{Remark}

\newcommand{\R}{\mathbf{R}}

\newcommand{\Z}{\mathbf{Z}}

\renewcommand{\tilde}{\widetilde}

\def\nxn{n\times n}

\def\sm#1{\{\pm1\}^{#1\times#1}}
\def\md{\mathrm{md}}

\def\tns{\otimes}
\def\va{{\bf a}}
\def\vb{{\bf b}}
\def\vc{{\bf c}}
\def\vd{{\bf d}}
\def\vx{{\bf x}}
\def\vy{{\bf y}}
\def\tp{{}^t\negthinspace}
\def\vo{{\mathbf \epsilon}}
\def\blk#1{{\mathbf{#1}}}
\def\Hom{\mathrm{Hom}(\R^n)}
\def\had#1{\mathrm{Had}\left(#1\right)}
\def\rs#1{\mathrm{RS}\left(#1\right)}

\def\ex#1{\mathrm{ex}\left(#1\right)}
\def\Nt{{N'}}
\def\tn{\Gamma(N)}
\def\rhomin{\rho_{\mathrm{min}}}
\def\slashfrac#1/#2{{{}^#1\negthinspace/\negthinspace_#2}}
\def\sf{\slashfrac}
\def\ph{\phantom}
\def\a{\alpha}
\def\s{\sigma}

\begin{document}

\parskip 6pt
\parindent 0pt
\baselineskip 14pt

\title{Large-determinant sign
matrices of order $\,\bf 4k+1\,$.}

\author[Orrick]{William P. Orrick}
\address{Department of Mathematics, Indiana University,
Bloomington, IN 47405}
\email{worrick@indiana.edu}

\author[Solomon]{Bruce Solomon}
\address{Department of Mathematics, Indiana University,
Bloomington, IN 47405}
\email{solomon@indiana.edu}

\date{First draft June 1, 2003. Last Typeset \today.}

\begin{abstract}
The Hadamard maximal determinant problem asks for the largest
$\,\nxn\,$ determinant with entries $\,\pm 1\,$. When
$\,n\equiv1\pmod 4\,$, the maximal excess construction of
Farmakis \& Kounias \cite{FK} has been the most successful
general method for constructing large (though seldom maximal)
determinants. For certain small $\,n\,$, however,
still larger determinants have been known; several new records
were recently reported in \cite{OSDS}. Here, we define
``3-normalized''  $\,\nxn\,$ Hadamard matrices, and construct
large-determinant matrices of order $\,n+1\,$ from them. Our
constructions account for most of the previous ``small $n$''
records, and set new records when
$\,n=37,\,49,\,65,\,73,\,77,\,85,\,93\,$, and $97$, most of
which are beyond the reach of the maximal excess technique. We
conjecture that our $\,n=37\,$ determinant,
$\,72\times9^{17}\!\times2^{36}\,$, achieves the global maximum.

\end{abstract}

\maketitle

\section{Introduction}\label{sec:intro}

\subsection*{Overview}
How big can the determinant of an $\,\nxn\,$ real matrix be,
given a bound on the size of its entries? This is the Hadamard
maximal determinant problem, and an easy argument reduces it to
that of maximizing the determinant on matrices with entries
$\,\pm1\,$.  From here on, we consider only such
matrices\footnote{Maximizing the determinant of an $\,\nxn\,$
$\,\pm1\,$ matrix is also equivalent to maximizing the
determinant over $\,\{0,1\}$-matrices of order $\,n-1\,$. See
\cite{Wi}}. Determinant-maximizing matrices of this type are
important in the design of experiments; they are called
a {D-optimal designs}.

When $\,n=4k+1\,$, the determinant of a $\,\pm1\,$ matrix can
never exceed the bound
\begin{equation}\label{eqn:barba}
B(n):=\left(n-1\right)^{n-1\over 2}\sqrt{2n-1}
\end{equation}
due to Barba \cite{Ba}. In these orders, the maximal excess
construction of Farmakis \& Kounias \cite{FK} has offered the
most successful general method for constructing determinants
achieving a large fraction of $\,B(n)\,$. For certain orders
such as $\,n=13,\,21,\,25,\,41,\,$ and $\,61\,$, however,
determinants beyond those attainable via maximal excess have
been known for some time \cite{R,CKM,Bh,BHH,Br}. In other
dimensions, specifically
$\,n=29,\,33,\,45,\,53,\,69,\,77,\,85\,$ and $\,93\,$, a
numerical ascent algorithm recently produced the largest known
determinants \cite{OSDS}. Here, we present a general
construction which accounts for many of the cases above. A
slight modification of our method accounts for the case
$\,n=29\,$. Except for certain maximizers which have been
obtained as incidence matrices for suitable block designs
($n=25\,,41$ and $61$), this leaves only the current $n=33$
record unexplained. In addition to clarifying the structure of
many current records, our construction improves the previous
records in dimensions $\,n=77,\,85\,$ and $93$, and sets new
records in dimensions $\,n=49,\,65,\,73\,$, and $97$. In a
sequel to the present paper, we describe the modified
construction which accounts for the case $\,n=29\,$, and show
that for $\,n=37\,$, it produces the determinant
$\,72\times9^{17}\!\times2^{36}=0.94B(37)\,$. We conjecture
that this value achieves the global maximum in its dimension.

We base our construction on \textit{3-normalized} Hadamard
matrices, which we define in \S\ref{sec:prelims} below. Starting
with a 3-normalized Hadamard matrix of order $\,4k\,$, we make
certain rank-1 modifications, and then construct a
large-determinant matrix of order $\,4k+1\,$ by adjoining an
additional row and column.

To optimize our method in a given dimension, we seek a
3-normalized Hadamard ``starting'' matrix whose excess, defined
as the sum of all matrix entries (see Definition \ref{def:ex}),
is as large as possible. In \S\ref{sec:bounds}, we derive an
upper bound on this excess. The determinants attainable by the
maximal excess construction are also bounded, given the bounds
on the excess of an Hadamard matrix first noted by Brown \&
Spencer \cite{BS}, and later improved by Kounias \& Farmakis
\cite{KF}. Comparison of the analogous bounds suggests that our
method has the potential to construct determinants larger than
any obtainable using the maximal excess method in arbitrarily
high dimensions.

Furthermore, for $\,n\,$ up to about 100, we find Hadamard
matrices that attain or closely approach our bound, thus
establishing the records mentioned above. On the other hand, we
do not yet know of any infinite family of Hadamard matrices that
attains our bound, while several such families are known for the
maximal excess construction. We should also note that
asymptotically, the largest fraction of Barba's bound that
either method can attain is $\,B(n)/\sqrt{2}\approx 0.71
B(n)\,$.  Our bound does exceed that of Kounias \& Farmakis by a
lower order term as $\,n\to \infty\,$, however, as shown in
Proposition \ref{prop:compare}.

\section{Preliminaries}\label{sec:prelims}

\subsection*{Maximal determinants}
Let $\,\sm{n}\,$ denote the space of all $\,\nxn\,$ \emph{sign
matrices}---matrices populated entirely by $\,\pm1$s.
This
space has cardinality $\,2^{n^2}<\infty\,$, so for
each $\,n\in|\Z|\,$, there exists a sign matrix of maximal
determinant.
Define the \emph{maximal
determinant function} accordingly:
\begin{equation}\label{eqn:md}
\md:|\Z|\to|\Z|\ ,\qquad
\md(n):=\max\left\{\det A\ :\ A \in\sm{n}\ \right\}\ .
\end{equation}
The exact value of $\,\md(n)\,$ has been established for all
$\,n\le 26\,$ except $\,n= 19,\,22\,$ and $23$. An elementary
argument first given by Hadamard shows that we always have
$\,\md(n)\le n^{n/2}\,$,  with equality iff the maximizing
matrix is \emph{Hadamard}, i.e. a sign matrix with mutually
orthogonal rows. A well-known conjecture states that Hadamard
matrices of order $\,n\,$ exist for all $\,n\equiv0\pmod4\,$, so
that the bound is achieved in these dimesions. This
conjecture is known to hold for all $\,n\equiv0\pmod4\,$ with
$\,n<428\,$, and many larger values as well~\cite{SY}.

The present paper focuses entirely on the case
$\,n\equiv1\pmod{4}\,$, for which $\,\md(n)\,$ remains largely
unknown. For $\,n\ge29\,$ with $\,n\equiv 1\pmod 4\,$, we only
know $\,\md(n)\,$ when $\,n\,$ takes the form $\,a^2+(a+1)^2\,$
when $\,a\,$ is an odd prime power \cite{Br}, and when $\,a=4\,$
\cite{BHH}. For these values of $\,n\,$, as well as $\,n=1\,$, 5,
13, and 25, the maximizing matrix is essentially the incidence
matrix of a suitable block design, and $\,\md(n)\,$ attains the
upper bound $\,B(n)\,$ for odd $\,n\,$ given in
(\ref{eqn:barba}) above. Conversely, equality cannot hold
here unless $\,n=a^2+(a+1)^2\,$ for some integer $\,a\,$, which
forces $\,n\equiv1\pmod 4\,$. Conjecturally, Barba's bound is
attained whenever $\,n\,$ has this form.

\subsection*{Tensor notation}
Given vectors $\,\vx,\,\vy\in\R^n\,$ with coordinates
$\,x_i,\,y_j\,$,  we will identify the tensor
product $\,\vx\tns\vy\,$ with the rank-1 linear transformation
whose action on $\,\R^n\,$ and matrix entries
$\,[\vx\tns\vy]_{ij}\,$ are characterized respectively by
\begin{equation}\label{eqn:tnsrDef}
(\vx\tns\vy)({\bf v}) := (\vy\cdot{\bf v})\,\vx\ ,
\qquad [\vx\tns\vy]_{ij} = x_i\,y_j\ .
\end{equation}
For any $\,A\in\Hom\,$, we then have
\begin{equation}\label{eqn:tnsrProp1}
A\circ(\vx\tns\vy) = (A\vx)\tns\vy\ ,\quad\text{and}\quad
(\vx\tns\vy)\circ A = \vx\tns(\tp A\vy)\ ,
\end{equation}
where $\,\tp A\,$ denotes the transpose of $\,A\,$.
In particular, when $\,A={\bf c}\tns{\bf d}\,$ for some $\,{\bf
c},{\bf d}\in\R^n\,$, we have
\begin{equation}\label{eqn:tnsProp2}
(\vx\tns\vy)\circ({\bf c}\tns{\bf d}) = (\vy\cdot{\bf
c})\vx\tns{\bf d}\ .
\end{equation}
Note also that
\begin{equation}\label{eqn:xpose}
\tp(\vx\tns\vy) = \vy\tns\vx\ .
\end{equation}
Using these facts, we get well-known formulae for the inverse and
determinant of the ``rank-1 update'' of a non-singular matrix:
\medskip

\begin{lem}[Sherman-Morrison formulae]\label{lem:rank1}
If $\,\vx,\vy\in\R^n\,$ and $\,A_{\nxn}\,$ is invertible,
then
\begin{eqnarray*}
\det(A+\vx\tns\vy) &=& (\det A) \,(1+ A^{-1}\vx\cdot\vy)\\
(A+\vx\tns\vy)^{-1} &=& A^{-1}- {A^{-1}\vx\tns (\tp A^{-1})\vy\over1+ A^{-1}\vx\cdot\vy}
\end{eqnarray*}
\end{lem}

\begin{proof}
First consider the special case $\,A=I_n\,$, and choose an
orthonormal basis $\,\{\tilde e_j\}\,$ for $\,\R^n\,$ with
$\,\tilde e_1 = \vx/|\vx|\,$. Then
\[
(\vx\tns\vy)(\tilde e_j) = (\vy\cdot\tilde e_j)\vx =
(\vy\cdot\tilde e_j)|\vx|\tilde e_1\ ,
\]
It follows that, relative to this basis, the matrix representing
$\,I_n+\vx\tns\vy\,$ has zeros everywhere except along the first
row, and the diagonal. In particular, it is upper triangular,
with all diagonal entries equal to 1 except the first, which is
$\,1+ \vx\cdot\vy\,$. The determinant of $\,I_n+\vx\tns\vy\,$
then reduces to this single entry, and we then have, for any
invertible $\,\nxn\,$ matrix $\,A\,$,
\begin{eqnarray*}
\det(A+\vx\tns\vy)
&=& (\det A)\,\det\left(I_n+A^{-1}(\vx\tns\vy)\right)\\
&=& (\det A)\,\det\left(I_n+(A^{-1}\vx)\tns\vy\right)\\
&=& (\det A)\,(1+A^{-1}\vx\cdot\vy)\ .
\end{eqnarray*}
This gives the stated formula for determinants.

One can verify the formula for the inverse directly: Multiply it
on the right by $\,(A_\vx\tns\vy)\,$, and simplify, using the
rules in \S \ref{sec:prelims}. The result, a straightforward
calculation, is $\,I_n\,$.
\end{proof}

\begin{definition}\label{def:on}
To streamline the display of vectors and matrices below, we
introduce the following convention: Given any integer $\,n\,$,
and any $\,k\in|\Z|\,$,
$\,\mathbf{n}_k:=n\,(1,1,\dots,1)\in\R^k\,$. Abusing notation
slightly, we also use $\,\blk{n}_i\,$ to indicate a sequence of
$\,k\,$ copies of the number $\,n\,$. For instance,
\[
(\,\blk{3}_{11},\,-\blk{1}_9\,)=
(\,\underbrace{3,\,3,\,\dots,\,3}_{11}\,,\,\underbrace{-1,\,-1,\,\dots,\,-1}_{9}\,)\in\R^{20}\
.
\]
\end{definition}

\subsection*{Excess}
If $\,A\in\sm{n}\,$, then the sum of the entries of $\,A\,$
measures the extent to which the $+1$s ``exceed'' the number of
$-1$s in $\,A\,$.

\begin{definition}[Excess and row sums]\label{def:ex}
The \emph{excess} of any vector or matrix $\,X\,$ is the sum of
its entries, and we denote it by $\,\ex{X}\,$. Note that for
a vector $\,\vx\in\R^k\,$ or an $\,n\times m\,$ matrix $\,A\,$
respectively, we have
\begin{equation}\label{eqn:ex}
\ex{\vx}:=\vx\cdot\blk{1}_k\ ,\qquad\ex{A}:=\blk{1}_n\cdot A\blk{1}_m\ .
\end{equation}
We call the excess of the $i$th row of a matrix its $i$th
\emph{row sum}.
\end{definition}

\begin{obs}\label{obs:ex}
If $\,\vx\in\R^n\,$ and $\,\vy\in\R^m\,$, then
\[
\ex{\vx\tns\vy} = \ex{\vx}\,\ex{\vy}\ .
\]
\end{obs}

This follows immediately from (\ref{eqn:tnsrDef}).

\subsection*{The injection $\,E:\R^{\nxn}\to\R^{(n+1)\times(n+1)}\,$}\label{ssec:E}\
Let $\,\R^{\nxn}\,$ denote the space of all real $\,\nxn\,$
matrices. Several of the facts above come together in relation
to the map $\,E:\R^{\nxn}\to\R^{(n+1)\times(n+1)}\,$, which adds
an additional first row and column to $\,A\,$ by prepending
$-1$ to each row of $\,A\,$, and then prepends an entire row of
$+1$s to that. In block form, we may express this as follows:
\begin{equation}\label{eqn:E}
E(A_{\nxn}) :=
\left(
\begin{array}{ll}\phantom{-}1 & \,\blk{1}_n\\
-\blk{1}_n & A
\end{array}
\right)\ .
\end{equation}

\begin{lem}\label{lem:detE}
For any square matrix $\,A\,$, we have
\[
\det\left(E(A)\right) = (\det A)\,\left(1+\ex{A^{-1}}\right)\ .
\]
\end{lem}
\begin{proof}
Adding the first row of $\,E(A)\,$ to each successive row
produces the matrix
\[
\left(
\begin{array}{lc}
1 & \blk{1}_n\\
\blk{0}_n & A+\blk{1}_n\tns\blk{1}_n
\end{array}
\right)\ ,
\]
without changing the determinant. The determinant of the latter
matrix is simply the determinant of its lower-right block,
however. We thus have
\[
\det\left(E(A)\right) = \det(A+\blk{1}_n\tns\blk{1}_n)\ .
\]
The desired result now follows immediately from Lemma
\ref{lem:rank1} and (\ref{eqn:ex}).
\end{proof}

\subsection*{Hadamard matrices.}
As mentioned above, a \emph{Hadamard matrix} is an
element $\,H\in\sm{n}\,$ whose rows are mutually orthogonal. We
denote the set of $\,\nxn\,$ Hadamard matrices by $\,\had{n}\,$.

Row and column permutations, and row and column negations,
clearly leave $\,\had{n}\,$ invariant. Hence

\begin{definition}\label{def:equiv}
Hadamard matrices $\,H_1,\,H_2\in\had{n}\,$ are \emph{equivalent}
if $\,H_1\,$ can be transformed into $\,H_2\,$ by permuting
and/or negating sets of rows and/or columns.
\end{definition}

We also note here that $\,\had{n}\subset\sm{n}\,$ is invariant
under transposition. Indeed, the defining condition for
$\,H\in\had{n}\,$ implies that $\,H\,\tp H=n\,I_n\,$, and
hence that $\,H^{-1} = \tp H/n\,$.

\begin{prop}\label{prop:stan}
If $\,n>2\,$ and $\,H\in\had{n}\,$, then $\,n=4k\,$ for some
$\,k\in\Z^+\,$, and using only column swaps and column
negations, we can put the first 3 rows of $\,H\,$ into
the following {standard} form
\begin{equation}\label{eqn:stan}
\left(
\begin{array}{cccc}
+&+&+&+\\
+&+&-&-\\
+&-&+&-
\end{array}
\right)\ .
\end{equation}
Here `$\,+$' and `$\,-$' respectively abbreviate $\,+\blk{1}_k\,$ and
$\,-\blk{1}_k\,$.  Moreover, if we denote the $i$th row sum of
$\,H\,$ by $\,r_i\,$, we have
\begin{equation}\label{eqn:quadratic}
\sum_{i=1}^{4k}r_i^2 = 16k^2\ .
\end{equation}
\end{prop}

\begin{proof}
To begin, note that by negating column $\,j\,$ of
$\,H=[h_{ij}]\,$ whenever $\,h_{1j}=-1\,$, we get an equivalent
matrix whose first row contains only $+1$'s. If we then sort the
columns of the resulting matrix in decreasing lexicographic
order, left to right, the first three rows of the resulting
$H$-equivalent matrix will necessarily take the form
\begin{equation*}
\left(
\begin{array}{cccc}
+\blk{1}_{k_1}&+\blk{1}_{k_2}&+\blk{1}_{k_3}&+\blk{1}_{k_4}\\
+\blk{1}_{k_1}&+\blk{1}_{k_2}&-\blk{1}_{k_3}&-\blk{1}_{k_4}\\
+\blk{1}_{k_1}&-\blk{1}_{k_2}&+\blk{1}_{k_3}&-\blk{1}_{k_4}
\end{array}
\right)\ ,
\end{equation*}
where $\,k_1+k_2+k_3+k_4=n\,$, and the mutual
orthogonality of these rows gives three additional
equations:
\begin{eqnarray*}
k_1+k_2-k_3-k_4&=&0\\
k_1-k_2+k_3-k_4&=&0\\
k_1-k_2-k_3+k_4&=&0
\end{eqnarray*}
The unique solution of the resulting $\,4\times 4\,$ system is
$\,k_1=k_2=k_3=k_4=n/4\,$, and since each $\,k_i\,$ is an
integer, we have $\,n/4\in\Z^+\,$. This proves the first
two claims.

To get the fact about row sums, note that the
the vector $\,(r_1,\,r_2,\,\dots,\,r_{4k})\,$ of row sums
coincides with $\,H\blk{1}_{4k}\,$. It follows that
\[
\sum_{i=1}^{4k}r_i^2 = H\blk{1}_{4k}\cdot
H\blk{1}_{4k}=\blk{1}_{4k}\cdot (\tp H H)\blk{1}_{4k} =
\blk{1}_{4k}\cdot (4k\blk{1}_{4k}) = (4k)^2\ .
\]
This completes the proof.
\end{proof}

Though the first three rows of any Hadamard matrix can always be
put in the form (\ref{eqn:stan}), a different
normalization will better suit our purposes.

\begin{definition}[3-normalization]\label{def:nor}
When the first three rows of a Hadamard matrix $\,H\,$
have the form
\begin{equation}\label{eqn:normform}
\left(
\begin{array}{cccc}
+&-&-&+\\
+&-&+&-\\
+&+&-&-
\end{array}
\right)\ ,
\end{equation}
\emph{and} all row sums of $\,H\,$ are non-negative, we say that
$\,H\,$ is \textit{3-normalized}.
\end{definition}
\goodbreak

\begin{prop}\label{prop:nor}
Every Hadamard matrix $\,H\in\had{n}\,$, $\,n>2\,$, can be
3-normalized using column swaps, column-negation, and row-negation
only. Moroever, every 3-normalized Hadamard matrix has the form
\smallskip
\begin{equation*}
N= \left(
\begin{array}{cccc}
+&-&-&+\\
+&-&+&-\\
+&+&-&-\\
\va_4&\vb_4&\vc_4&\vd_4\\
\vdots&\vdots&\vdots&\vdots\\
\va_n&\vb_n&\vc_n&\vd_n\\
\end{array}
\right)\ .
\end{equation*}
\smallskip
Here `$\,+$' and `$\,-$' respectively abbreviate
$\,\pm\blk{1}_k\,$, and
$\,\va_i,\,\vb_i,\,\vc_i,\,\vd_i\in\{\pm1\}^{n/4}\,$, and
satisfy, for each $\,i=4,\dots,n\,$,
\begin{equation}\label{eqn:rsums}
\ex{\va_i}=\ex{\vb_i}=\ex{\vc_i}=\ex{\vd_i}\ .
\end{equation}
For the row sums $\,r_i\,$, we have $\,r_1=r_2=r_3=0\,$, and, for
$\,i>3\,$, $\,r_i\equiv n\pmod8\,$.
If $\,n\equiv 4\pmod
8\,$, $\,r_i\le n-8\,$ for all rows, unless $\,n=4\,$.

Finally, we have $\,\ex{N}\equiv n\pmod 8\,$, and when $\,n\equiv
0\pmod 8\,$, $\,\ex{N}\equiv n\pmod{16}\,$.
\end{prop}

\begin{proof}
By Proposition \ref{prop:stan}, we can write $\,n=4k\,$ and put
the first three rows of $\,H\,$ into the form (\ref{eqn:stan})
using only column swaps and negations. By further negating
columns $\,k+1\,$ through $\,3k\,$, we obtain an $H$-equivalent
matrix whose first three rows have the form
(\ref{eqn:normform}). Negating all rows having negative row
sums, we now get a 3-normalized Hadamard matrix equivalent to
$H$. This proves the Lemma's first statement.

To get (\ref{eqn:rsums}), define $\,a_i=\ex{\va_i}\,$,
$\,b_i=\ex{\vb_i}\,$, etc., and recall that for any fixed
$\,i=4,5,\dots,n\,$, row $\,i\,$ is orthogonal
to each of the first three rows. Consequently,
\begin{eqnarray*}
a_i-b_i-c_i+d_i&=&0\\
a_i-b_i+c_i-d_i&=&0\\
a_i+b_i-c_i-d_i&=&0\ .
\end{eqnarray*}
The general solution here is $\,a_i=b_i=c_i=d_i\,$, which proves
(\ref{eqn:rsums}). It also shows that $\,r_i\equiv0\pmod4\,$ for
all $\,i\,$.  In fact, the first three row sums clearly equal 0,
and to get the more precise claim about rows 4 through $\,n\,$,
simply note that since $\,\va_i\,$ has $\,k\,$ entries, each
congruent to 1 mod 2, $\,\ex{\va_i}\,$ has the same parity as
$\,k\,$ does. It follows that $\,r_i=4\,\ex{\va_i}\equiv 0\,$ or
$\,4\pmod8\,$, depending on whether $\,k\,$ is even or odd
respectively.

For odd $\,k>1\,$, no row can consist entirely of 1s (giving row
sum $\,n\,$), because the mutual orthogonality of rows would
then force all remaining row sums to equal
$\,0\not\equiv4\pmod8\,$. Therefore, no row sum exceeds
$\,n-8\,$ in this case, as claimed.

Since we have an odd number ($n-3$) of non-zero row sums, all
congruent to $\,n\pmod 8\,$, it follows that $\,\ex{N}\equiv
n\pmod 8\,$.

Finally, since $\,r_i\equiv n\pmod 8\,$, we may write
$\,r_i=8l_i\,$ and $\,n=8m\,$ when $\,n\equiv 0\pmod 8\,$. By
definition of the the excess, and (\ref{eqn:quadratic}) above,
we then have
\begin{equation*}
{\ex{N}\over 8} = \sum_{i=1}^n l_i\qquad\text{and}\qquad
m^2 = \sum_{i=1}^n l_i^2\ .
\end{equation*}
Since $\,x^2\equiv x\pmod2\,$ for any integer $\,x\,$, the two
sums above, and $\,m\,$, all have the same parity. We conclude
that $\,\ex{N}\equiv 8m\pmod{16}\,$, as claimed.

\end{proof}

\begin{remark}
When $\,n\ge 20\,$ and $\,n\equiv 4\pmod 8\,$, the proof of Lemma
2.2 in \cite{Ha} generalizes to show that the largest possible
row sum allowed by Proposition \ref{prop:nor} above, namely
$\,n-8\,$, can be achieved by at most one row. Modulo column
negations and permutations, such a row forms, along with the
three initial rows of the 3-normalized matrix to which it
belongs, a \textit{Hall set}, as defined by Kimura \cite{Ki}.
\end{remark}

\section{A construction}\label{sec:construction}

\subsection*{Constructing $\,\tn\,$}
We now describe the basic construction we use to produce
large-determinant sign matrices of order $\,n+1\,$.

\begin{description}

\item[Step 1]
Select a \emph{3-normalized} Hadamard matrix $\,N\in\had{n}\,$,
$\,n=4k\,$.
\medskip
\item[Step 2]
Alter $\,N\,$ by negating its first $\,k\,$ {columns} to
produce a new Hadamard matrix $\,\Nt\,$. By Lemma \ref{prop:nor}
the first three rows of $\,\Nt \,$ will take the form
\begin{equation}
\left(
\begin{array}{cccc}
-&-&-&+\\
-&-&+&-\\
-&+&-&-
\end{array}
\right)\ .
\end{equation}
\medskip

\item[Step 3]
Make the rank-1 modification
\begin{equation*}
\Nt \longrightarrow \Nt  + 2\,\vo_3\tns\vo_{k}=:N''\ ,
\end{equation*}
where, for any $\,m\in|\Z|\,$ we define $\,\vo_m :=
(\blk{1}_m,\,\blk{0}_{n-m})\in\R^n\,$. This modifies $\,\Nt\,$
by changing the $-1$s in its upper left $\,3\times k\,$ block to
$+1$s. Note also that $\,N''\,$ is no longer Hadamard.
\medskip

\item[Step 4]
Construct the final sign matrix $\,\tn \,$ of size $\,n+1\,$ by
applying the injection $\,E\,$ defined in \S\ref{sec:prelims}:
\[
\tn=E(N'')\ .
\]
\end{description}
\medskip

The following result suggests the efficacy of this construction.

\begin{thm}\label{thm:construction}
Let $\,N\,$ be a {\em 3-normalized} Hadamard matrix of order
$\,n=4k\,$. Then the sign matrix
$\,\tn\,$ of order $\,n+1\,$ described above has determinant
\[
\det\tn =  (\det N)\left(2+{\ex{N}\over
n}\right)=\pm n^{n/2}\left(2+{\ex{N}\over n}\right)\ .
\]
\end{thm}

\begin{proof}
Using the notation of our construction above, we can express
$\,\det \tn \,$ in terms of $\,N''\,$ by applying Lemma
\ref{lem:detE} to get:
\begin{equation}\label{eqn:detR1}
\det \tn =\det N'' \left(1+\ex{N''^{-1}}\right)\ .
\end{equation}
We proceed to expand the right-hand side in terms of $\,N\,$
itself.

First, given the definition of $\,N''\,$ in Step 3 of our
construction, Lemma \ref{lem:rank1} enables us to write
\begin{equation}\label{eqn:detN''}
\det N''  =
(\det\Nt)\left(1+2\,\Nt^{-1}\vo_3\cdot\vo_{k}\right)\ .
\end{equation}
Since $\,\Nt\,$ is Hadamard, we have $\,\Nt^{-1}=\tp
\Nt/n\,$, and hence $\,\Nt^{-1}\vo_3= \tp \Nt\vo_3/n\,$. But
$\,\tp \Nt\vo_3\,$ simply sums the the first three columns of
$\,\tp \Nt\,$, ie. the first three \emph{rows} of $\,\Nt\,$
itself. Step 2 of our construction displays these rows
explicitly, and we easily deduce
\begin{equation}\label{eqn:Nt3}
\Nt^{-1}\vo_3 = {1\over n}
\left(\,-\mathbf{3}_{k},\,-\blk{1}_{3k}\,\right)\ .
\end{equation}
It follows that
\begin{equation}\label{eqn:denom}
1+ 2\,\Nt^{-1}\vo_3\cdot\vo_{k} = 1 + {2\over
n}\,
\left(
\begin{array}{l}
-\blk{3}_{k}\\
-\blk{1}_{3k}
\end{array}
\right)
\cdot
\left(
\begin{array}{l}
\blk{1}_{k}\\
\blk{0}_{3k}
\end{array}
\right) = 1- {2\over
n}\cdot{3k} = -{1\over 2}\ ,
\end{equation}
and consequently,
\begin{equation}\label{eqn:detM2}
\det N'' = (-1)^{k+1}\, {\det N\over 2}\ .
\end{equation}
The factor $\,(-1)^k\,$ appears because we have replaced
$\,\Nt\,$ by $\,N\,$, a matter of negating $\,k\,$ columns.

We next calculate $\,\ex{N''^{-1}}\,$.

To begin, invert $\,N''=\Nt+2\,\vo_3\,\tns\,\vo_k\,$ using Lemma
\ref{lem:rank1}, and compute
\begin{eqnarray}\label{eqn:exMinv}
\ex{N''^{-1}}
&=& \ex{\Nt^{-1}}- 2\,{\ex{\Nt^{-1}\vo_{3}\,\tns\,(\tp\Nt)^{-1}\vo_k}\over
1+2\,\Nt^{-1}\vo_3\cdot\vo_k}\nonumber\\
&=& \ex{\Nt^{-1}}+ 4\,{\ex{\Nt^{-1}\vo_{3}}\cdot\ex{(\tp\Nt)^{-1}\vo_k}}\ .
\end{eqnarray}
Here we used (\ref{eqn:denom}) to evaluate the
denominator, and Observation \ref{obs:ex} to factor the excess of
the tensor product.

Now unpack the excess terms on the right above. For
$\,i=1,\,2,\,\dots,\,n\,$, let $\,r_i\,$ and $\,r'_i\,$ denote
the $i$th row sums of $\,N\,$ and $\,\Nt\,$ respectively. Since
$\,\Nt\,$ is Hadamard, we may then write
\begin{equation*}
\ex{\Nt^{-1}} = {1\over n}\,\ex{\Nt} = {1\over
n}\,\sum_{i=1}^n\,r'_i\ .
\end{equation*}
We obtained $\,\Nt\,$ from a 3-normalized Hadamard matrix $\,N\,$
by negating columns 1 through $\,k\,$. So by the definition of
3-normalized form, we have $\,r'_i = -2k = -n/2\,$ when
$\,i=1,2,3\,$. When $\,i\ge 4\,$ Proposition \ref{prop:nor}
similarly implies $\,r'_i = r_i/2\,$. Finally, since
$\,r_1=r_2=r_3=0\,$, we have $\,\sum_{i=4}^n r_i =
\sum_{i=1}^n r_i = \ex{N}\,$. Thus
\begin{equation}\label{eqn:exNtinv}
\ex{\Nt^{-1}} = {1\over
n}\cdot\left(\sum_{i=1}^3r'_i+\sum_{i=4}^nr'_i\right)= {1\over
n}\cdot\left({-3n+\ex{N}\over 2}\right)= {\ex{N}\over
2n}-{3\over 2}\ .
\end{equation}

Given these facts, we now easily analyze $\,\ex{\Nt^{-1}\vo_3}\,$
using (\ref{eqn:Nt3}):
\begin{equation}\label{eqn:exNt3}
\ex{\Nt^{-1}\vo_3}={1\over n}\,\left((-3)\,{n\over 4} +
(-1){3n\over 4}\,\right) = -{6\over 4}=-{3\over 2}\ ,
\end{equation}

Lastly, we compute $\,\ex{(\tp\Nt)^{-1}\vo_k}\,$, recalling
again that we produce $\,\Nt\,$ from $\,N\in\had{n}\,$ by
negating columns 1 through $\,k\,$, so that
\begin{eqnarray}\label{eqn:sumOcols}
(\tp\Nt)^{-1}\vo_k = {1\over n}\Nt\vo_k
&=&\phantom{-}{1\over n}\cdot(\text{sum of the first $\,k\,$ columns of $\,\Nt\,$})\nonumber\\
&=&-{1\over n}\cdot(\text{sum of the first $\,k\,$ columns of $\,N\,$})\nonumber\ .
\end{eqnarray}
Since $\,N\,$ is 3-normalized, rows 1, 2, and 3 have only $+1$s
in these columns. The first three rows therefore contribute
$\,k=n/4\,$ each to the sum of columns above.
Proposition~\ref{prop:nor} further guarantees that rows 4
through $\,n\,$ each contribute exactly 1/4th of their full row
sum. Denoting the latter row sums by $\,r_i\,$ as above, we may
therefore write
\begin{equation*}
(\tp\Nt)^{-1}\vo_k = -{1\over
4n}\,\left(\blk{n}_{3},\,r_4,\,r_5,\,\dots,\,r_n\right)\ ,
\end{equation*}
whence
\begin{equation}\label{eqn:exNtn/4}
\ex{(\tp\Nt)^{-1}\vo_k}=-{1\over
4n}\,\left(3\,n+\ex{N}\right)=-{\ex{N}\over 4n}-{3\over 4}\ ,
\end{equation}
since $\,\sum_{i=4}^nr_i = \ex{N}\,$, as noted in deriving
(\ref{eqn:exNtinv}) above.

Finally, we assemble the various facts above to calculate
$\,\det\tn\,$. First, combine equations (\ref{eqn:detR1}),
(\ref{eqn:detM2}) and (\ref{eqn:exMinv}), to get
\begin{equation*}
\det\tn =(-1)^{k+1}
(\det N)\cdot{1\over2}\cdot\left(1+\ex{\Nt^{-1}} +
4\,\ex{\Nt^{-1}\vo_3}\cdot\ex{(\tp\Nt)^{-1}\vo_k}\,\right)\ .
\end{equation*}
Then replace the excess terms here with the results of
equations (\ref{eqn:exNtinv}), (\ref{eqn:exNt3}), and
(\ref{eqn:exNtn/4}). Routine simplification then verifies the
theorem.
\end{proof}

\section{Upper bounds}\label{sec:bounds}

Compare the formula in Theorem~\ref{thm:construction} above with
the analogous formula we get by applying Lemma \ref{lem:detE} to
an arbitrary $\,\nxn\,$ Hadamard matrix $\,H\,$:
\begin{equation}\label{eqn:fkdet}
\det\left(E(H)\right) =
(\det H)\left(1+\frac{1}{n}\,\ex{H}\right)\ .
\end{equation}
The latter formula clearly shows that we maximize $\,\det E(H)\,$
over all $\,H\in\had{n}\,$ by maximizing $\,\ex{H}\,$. Indeed,
this constitutes the maximal excess technique of Farmakis \&
Kounias \cite{FK} we mention in our introduction.

By comparison, Theorem~\ref{thm:construction} above shows that
we can improve on that technique in any dimension
$\,n\,$ where there exists a \emph{3-normalized} $\,\nxn\,$
Hadamard matrix $\,N\,$ such that
\begin{equation}\label{eqn:compare}
\ex{N}+n>\ex{H}\quad\text{for all $\,H\in\had{n}\,$}\ .
\end{equation}
Of course, the 3-normalized matrices form a relatively small
subset in any equivalence class of Hadamard matrices, and one
expects the excess function to reach a smaller max on this
subset than on the entire class. So it is not obvious that we
can attain the desirable situation expressed by
(\ref{eqn:compare}).

We have found, however, that 3-normalized Hadamard matrices
satisfying (\ref{eqn:compare}) do exist, at least in dimensions
$\,n< 100\,$. The excesses of arbitrary Hadamard matrices and
3-normalized Hadamard matrices satisfy certain natural upper
bounds. We discuss both bounds below, and compare them. We find
that for large $\,n\,$, the bounds themselves \textit{do}
satisfy (\ref{eqn:compare}), offering the hope that one could
find an actual 3-normalized matrix that does so.

In this section, we briefly review known bounds on the excess of
a general Hadamard matrix $\,H\in\had{n}\,$, then derive a new
upper bound for the excess of a 3-normalized Hadamard matrix.
Our bound clarifies what we can (and cannot) hope to obtain from
the construction in \S \ref{sec:construction}, and indicates how
to optimize the search for suitable 3-normalized Hadamard
matrices, which we take up in \S \ref{sec:examples}.

Proposition \ref{prop:compare} at the end of this section gives
our comparison of the two bounds.

\subsection*{The $\,n^{3/2}\,$ bound on maximal excess}
The excess of a general Hadamard matrix $\,H\in\had{n}\,$ can
never exceed $\,n^{3/2}\,$. Indeed, the Cauchy-Schwartz
inequality gives
\[
\ex{H} = \blk{1}_n\cdot H\blk{1}_n\le
|\blk{1}_n|\cdot |H\blk{1}_n|=\sqrt{n}\cdot\sqrt{n}\cdot\sqrt{n}\
,
\]
since $\,H/\sqrt{n}\,$ is orthogonal. This fact was apparently
first observed by Brown \& Spencer \cite{BS}, and independently,
by Best \cite{B}. The bound is sharp in the sense that regular
Hadamard matrices actually attain this excess value, and such
matrices occur in infinitely many orders (for instance, order
$\,n^2\,$ whenever $\,\had{n}\,$ is non-empty \cite{GS}).

In \cite{KF}, Kounias \& Farmakis derive an upper bound smaller
than $\,n^{3/2}\,$ when $\,n\,$ is not a perfect square, and
tabulate their bound versus the largest known excess for all
orders up to $\,n=100\,$ \cite{FK}. For convenience, we
reproduce their list of bounds in Table 1, \S
\ref{sec:examples}.

\subsection*{A 3-normalized excess bound}
If $\,N\,$ is a 3-normalized Hadamard matrix of order $\,n>2\,$,
Prop.~\ref{prop:nor} shows that the first three rows of
$\,N\,$ have row sum zero, while the remaining row sums belong
to the set
\begin{equation*}\label{eqn:rsumsets}
\rs{n} = \left\{
\begin{array}{ll}
\{0,8,16,\dots,n\}\ ,&\text{if }\ n\equiv0\pmod8\ ,\\
\{4\,\}\ ,&\text{if }\ n=4\,,\\
\{4,12,20,\dots,n-8\}\ ,&\text{if }\ n\equiv4\pmod8\ , n>4\ .
\end{array}\right.
\end{equation*}
\medskip

\begin{obs}\label{obs:sums}
Consider the last $\,n-3\,$ rows of a {3-normalized}
$\,N\in\had{n}\,$. For each $\,r\in\rs{n}\,$, let $\,n_r\,$
count how many of these rows have row sum $\,r\,$. Then
\begin{equation*}
\sum_{\rs{n}}\,n_r =n-3\ ,\quad
\sum_{\rs{n}}r\,n_r=\ex{N}\ ,\quad\text{and}\quad
\sum_{\rs{n}}r^2\,n_r= n^2\ .
\end{equation*}
\end{obs}

\begin{proof}
The first identity is trivial. The second totals the row
sums---the excesses---of the last $\,n-3\,$ rows. Since the
first three rows of $\,N\,$ contribute zero excess, this gives
$\,\ex{N}\,$. Equation (\ref{eqn:quadratic}) from
Prop.~\ref{prop:stan} evaluates the third sum; we record it in
slightly different notation here for convenience only.
\end{proof}

\subsection*{Notation}
In the theorem below, the notations
\[
\lceil x \rceil\quad\text{and}\quad\lceil x\rceil_+
\]
respectively denote the \emph{ceiling} and \emph{positive}
ceiling of $\,x\,$, ie. least integer, and the least
\emph{positive} integer, respectively, which exceeds or
equals $\,x\,$.
\medskip

\begin{thm}\label{thm:bounds}
Suppose $\,N\in\had{n}\,$ is \emph{3-normalized}, with $\,n\ge4\,$.
Then $\,\ex{N}\le \nu_n^*\,$, where
\begin{equation}\label{eqn:key}
\nu_n^*:= {\bar\rho_n(n-3)\over 2}+{(n-4)(n-12)\over
2\bar\rho_n}\ ,
\end{equation}
and
\begin{equation*}
\bar\rho_n:=
\left\{
\begin{array}{ll}
  8\left\lceil\displaystyle{n\over 8\sqrt{n-3}}\right\rceil-4\ ,&\mbox{if }\ n\equiv0\pmod8\ ,\\
  &\\
  8\left\lceil\displaystyle{n\over
8\sqrt{n-3}}-{1\over2}\right\rceil_+\ , &\mbox{if }\
n\equiv4\pmod8\ .
\end{array}\right.
\end{equation*}
Equality holds in (\ref{eqn:key}) iff the last $\,n-3\,$ row sums
of $\,N\,$ take the values $\,\bar\rho_n\pm 4\,$ {only}.
\end{thm}

\begin{proof}
Our main argument degenerates for the cases $\,n=4,\,8$ and $12$;
we consequently argue the latter separately after handling the
generic case $\,n>12\,$.

\subsection*{The case $\mathbf{n>12}$}
As in Observation \ref{obs:sums}, we let $\,n_r\,$ count how
many, among the last $\,n-3\,$ rows of $\,N\,$, have row sum
$\,r\,$ for each $\,r\in\rs{n}\,$. We then adapt a technique
appearing in Kounias and Farmakis \cite{KF} by considering the
following function $\,F\,$ of $\,\rho\in\rs{n}\,$:
\begin{equation}\label{eqn:trick}
F(\rho):=\sum_{\rs{n}} n_r (\rho-r)(\rho+8-r)\ge 0\ .
\end{equation}
The inequality here is crucial. It holds because the
$\,n_r$ are all non-negative, while $\,(\rho-r)\,$ and
$\,(\rho+8-r)\,$ necessarily give consecutive integer multiples
of 8, and hence a non-negative product.

On the other hand, $\,F(\rho)\,$ expands as
\begin{equation*}
F(\rho)=\rho(\rho+8)\sum n_r-(2\rho+8)\sum r\,n_r +\sum r^2 n_r\
,
\end{equation*}
and we can evaluate all three sums explicitly using
Observation \ref{obs:sums}. After careful simplification, we
may then rewrite the estimate (\ref{eqn:trick}) above as
\begin{equation}\label{eqn:gn}
\ex{N}\le G_n(\rho):={(\rho+4)(n-3)\over 2}
+{(n-4)(n-12)\over2(\rho+4)}\ .
\end{equation}
Note too, for later use, that equality obtains here iff
$\,F(\rho)=0\,$.

To extract the strongest possible result from (\ref{eqn:gn}), we
now minimize $\,G_n(\rho)\,$ over $\,\rs{n}\,$. Here, the
assumption $\,n>12\,$ comes into play; one easily sees that in
this case $\,G_n\,$ is concave up, hence attains its minimum on
$\,\rs{n}\,$ at the smallest $\,\rho\in\rs{n}\,$ for which the
forward difference $\,G_n(\rho+8)-G_n(\rho)\,$ is non-negative.
Accordingly, we compute
\begin{equation*}
G_n(\rho+8)-G_n(\rho) =
4\,{(\rho+8)^2(n-3)-n^2\over(\rho+4)(\rho+12)}\ ,
\end{equation*}
and observe that the denominator here is positive for all
$\,\rho\in\rs{n}\,$. This makes non-negativity of the forward
difference on $\,\rs{n}\,$ equivalent to
\begin{equation*}
\rho\ge {n\over \sqrt{n-3}}-8\ .
\end{equation*}
Since $\,\rs{n}\,$ depends on the residue class of $\,n\,$ mod 8,
so does the minimal solution $\,\rhomin\in\rs{n}\,$ for this
inequality, but a little thought shows that we can express it in
the following way:
\begin{equation*}
\rhomin =\left\{
\begin{array}{ll}
8\left\lceil\displaystyle{n\over8\sqrt{n-3}}\right\rceil-8\ ,&\text{if }\
n\equiv0\pmod8\\
&\\
8\left\lceil\displaystyle{n\over8\sqrt{n-3}}-{1\over2}\right\rceil-4\
,&\text{if }\ n\equiv4\pmod8\ .
\end{array}\right.
\end{equation*}

By evaluating (\ref{eqn:gn}) with $\,\rho=\rhomin\,$ and then
substituting $\,\bar\rho_n:=\rhomin+4\,$, we now obtain
(\ref{eqn:key}), the main result of our theorem.

It remains to analyze the condition for equality there. According
to the remark immediately following (\ref{eqn:gn}), this
condition corresponds to $\,F(\rhomin)=0\,$, and by
(\ref{eqn:trick}), the latter will happen if and only if we
have
\[
n_r(\rhomin-r)(\rhomin+8-r)=0\quad\text{for all}\ r\in\rs{n}\ .
\]
The factor $\,(\rhomin-r)(\rhomin+8-r)\,$ here vanishes when and
only when $\,r=\rhomin\,$ or $\,r=\rhomin+8\,$, or equivalently,
when $\,r=\bar\rho_n\pm 4\,$. Equality in (\ref{eqn:key})
consequently forces $\,n_r=0\,$ for all other $\,r\in\rs{n}\,$,
exactly as claimed.

\subsection*{The cases $\mathbf{n=4,8}$ and $\mathbf{12}$}
Since $\,\rs{4}=\rs{12}=\{4\}\,$, the excess of a 3-normalized
Hadamard matrix $\,N\,$ will equal $(4-3)\times 4=4$ if
$\,N\in\had{4}\,$, and $(12-3)\times 4=36$ if
$\,N\in\had{12}\,$. The reader will easily check that our
theorem accurately predicts these values and appropriately
assigns $\,\bar\rho_n=8\,$.

When $\,n=8\,$, we have $\,\rs{8}=\{0,8\}\,$. In fact, a
3-normalized matrix $\,N\in\had{8}\,$ will have exactly one row
with sum 8 and seven with row sum 0. For, on the one hand,
all eight rows cannot sum to zero; this would make them
linearly dependent. On the other hand, all entries in a row
with sum 8 must be $+1$s, so even two such rows create a
dependency. We consequently have $\,\ex{N}=8\,$, again
confirming the theorem.
\end{proof}

\begin{remark}\label{rem:n=80}
According to Prop.~\ref{prop:nor}, the excess of a 3-normalized
matrix $\,N\in\had{n}\,$ always has certain congruence
properties. The bound $\,\nu_n^*\,$ in Theorem \ref{thm:bounds},
however, may not have these properties, and in such cases, we
cannot achieve  $\,\ex{N}=\nu_n^*\,$. This occurs for the first
time, for instance when $\,n=80\equiv0\pmod{16}\,$. Then
$\,\bar\rho_n=12\,$ and $\,\nu_n^*=677\sf1/3$. By Proposition
\ref{prop:nor}, however, the excess of any actual 3-normalized
$\,N\in\had{80}\,$ must be divisible by 16. Since the largest
multiple of $\,16\,$ below $\,\nu_{80}^*\,$ is $\,672\,$, the
latter provides a sharper upper bound for $\,\ex{N}\,$ when
$\,n=80\,$. We know no 3-normalized $\,N\in\had{80}\,$ with
excess above 624, but if there exists a matrix
$\,N\in\had{80}\,$ with  $\,n_0=1\,$ (ie. 1 row with excess 0),
$\,n_8=68\,$ and $\,n_{16}=8\,$, it has excess 672, thus
attaining the sharper bound.

In general, when this situation arises (the next such occurs when
$\,n=104\,$), our bound can be sharpened similarly. The precise
formulation, using $\,\lfloor x\rfloor\,$ to denote the integer
part of $\,x\,$ is as follows:
\begin{equation*}
\ex{N}\le
\begin{cases}
8\left\lfloor\displaystyle{{\nu_n^*\over 8} - \frac{1}{2}}\right\rfloor+4 & \text{if $\,n\equiv4\pmod 8\,$,}\\
&\\
16\left\lfloor\displaystyle{{\nu_n^*\over 16} - \frac{1}{2}}\right\rfloor+8 & \text{if $\,n\equiv8\pmod{16}\,$,}\\
&\\
16\left\lfloor\displaystyle{\nu_n^*\over 16}\right\rfloor &
\text{if $\,n\equiv0\pmod{16}\,$.}
\end{cases}
\end{equation*}
We should emphasize, however, that even this sharpened bound may
not be attained for all $\,n\,$. Indeed, it remains unknown
whether $\,\had{n}\,$ is non-empty for all $\,n\equiv 0\pmod
4\,$.
\end{remark}
 
\subsection*{Comparison of bounds}
As discussed at the beginning of this section, our construction
(Theorem \ref{thm:construction}) produces a determinant larger
than any the maximal excess method can construct whenever there
exists a 3-normalized Hadamard matrix $\,N\in\had{n}\,$ for
which (\ref{eqn:compare}) holds. Given the $\,n^{3/2}\,$ bound
on maximal excess, (\ref{eqn:compare}) certainly holds whenever
\begin{equation}\label{eqn:compare2}
\ex{N}+n> n^{3/2}\ .
\end{equation}
We now show that Theorem \ref{thm:bounds} presents no obstruction
to the existence of such a matrix $\,N\,$, except in the very
lowest dimensions. Recall that $\,\nu_n^*\,$, defined in that
theorem, denotes our upper bound for the excess of a 3-normalized
Hadamard matrix in $\,\had{n}\,$.

\begin{prop}\label{prop:compare}
For all $\,n\ge12\,$ we have
\begin{equation}\label{eqn:maxmin}
\sqrt{(n-3)(n-4)(n-12)}\le \nu_n^* \le n\sqrt{n-3}\ ,
\end{equation}
so that $\,\nu_n^* = n^{3/2} - O(\sqrt{n}\,)\,$, and
$\,(\nu_n^*+n)-n^{3/2}\sim n\,$ as $\,n\to \infty\,$ in the
sense that
\begin{equation}\label{eqn:asymp}
\lim_{n\to\infty} {(\nu_n^*+n) - n^{3/2}\over n} = 1\ .
\end{equation}
In fact, for all $\,n\ge 88\,$, we have
\begin{equation}\label{eqn:88}
\nu^*_n+n> n^{3/2}\ ,
\end{equation}
\end{prop}

\begin{proof}
It is clear from the definition of $\,\nu_n^*\,$ (and that fact
that $\,\bar\rho_n>0\,$) that we can bound $\,\nu_n^*\,$ from
below by the minimum of the function
\begin{equation}\label{eqn:nu}
G_n(x)={(n-3)x\over 2} + {(n-4)(n-12)\over 2x}\ ,\qquad x\in
(0,\infty)\ .
\end{equation}
(We minimized essentially the same function in proving Theorem
\ref{thm:bounds}, but over a discrete domain.)

When $\,n>12\,$, the function is concave up, with one critical
point on $\,(0,\infty)\,$.  Computing the minimum explicitly on
this interval, one easily obtains the lower bound in
(\ref{eqn:maxmin}) when $\,n>12\,$. The bound is trivial when
$\,n=12\,$.

To get the upper bound $\,\nu_n^*<n\sqrt{n-3}\,$, we note that
the definition of $\,\bar\rho_n\,$ also implies
\[
{n\over \sqrt{n-3}}+4\ge \bar\rho_n\ge {n\over \sqrt{n-3}}-4\ .
\]
If we restrict $\,x\,$ to this interval, the maximum of
$\,G_n(x)\,$ must occur at an endpoint, and an
easy calculation shows that
$\,G_n(x)=n\sqrt{n-3}\,$ at both endpoints. This value
consequently bounds $\,\nu_n^*\,$ from above, as claimed.

The asymptotic result $\,\nu_n^*=n^{3/2}-O(\sqrt{n})\,$ now
follows. In fact, one easily uses the the bounds in
(\ref{eqn:maxmin}) to compute
\[
-{3\over2}\ge\limsup_{n\to+\infty}\,{\nu_n^*-n^{3/2}\over
\sqrt{n}}\ge \liminf_{n\to+\infty}\,{\nu_n^*-n^{3/2}\over
\sqrt{n}} \ge -{19\over 2}\ .
\]
These facts make (\ref{eqn:asymp}) obvious, and imply
that (\ref{eqn:88}) holds for sufficiently large $\,n\,$.

To show more precisely that (\ref{eqn:88}) holds whenever
$\,n\ge 88\,$, it will suffice to do so with $\,\nu_n^*\,$
replaced by the lower bound in (\ref{eqn:maxmin}); that is to
show
\[
\sqrt{(n-3)(n-4)(n-12)} >n^{3/2}-n\quad\text{for all
$\,n\ge88\,$}\ .
\]
Squaring both sides and simplifying the result we find that this
inequality is equivalent to
\[
(n^{5/2}-10n^2) + (48n-72)>0\ .
\]
When $\,\sqrt{n}> 10\,$, the first term in parentheses here is
clearly non-negative, and the second is positive. We
therefore see that (\ref{eqn:88}) holds for all $\,n\ge 100\,$.
By hand, one checks that it holds for $\,n=96\,,92\,$ and
$\,88\,$ as well, but not for $\,n=84\,$.
\end{proof}

As mentioned earlier, Kounias \& Farmakis bound the maximal
excess of an order-$n$ Hadamard by a number
$\,\sigma_n^*<n^{3/2}\,$ when $\,n\,$ is not a perfect square.
For $\,n<88\,$, our bound $\,\nu_n^*+n\,$ dips below
$\,n^{3/2}\,$, but it still exceeds $\,\sigma_n^*\,$ in most
cases. Indeed, for all $\,n\equiv 4\pmod 8\,$, we have
$\,\nu_n^*+n\ge \sigma_n^*\,$, with equality only for perfect
square $\,n\,$. For $\,n\equiv0\pmod 8\,$, we have
$\,\nu_n^*+n>\sigma_n^*\,$ for all $\,n\ge 44\,$. For
$\,n\le100\,$ we display the exact values of $\,\sigma_n^*\,$
and $\,\nu_n^*+n\,$ in Table 1 below.

\section{Examples}\label{sec:examples}

Theorem \ref{thm:bounds} bounds the excess of a 3-normalized
Hadmard matrix from above, and by virtue of Theorem
\ref{thm:construction} this provides an upper bound for the
determinants we can hope to achieve using the construction
described at the beginning of \S \ref{sec:construction}.
These facts would hold little interest if the construction did
not actually produce examples that attain, or at least approach
these upper bounds, thereby providing good lower bounds
for the maximal determinant problem in the corresponding
dimensions.

Our construction does produce numerous examples of this type.
With regard to low orders, it generates globally
determinant-maximizing sign matrices of orders $\,n=5$, $13$, and
$21$. By way of comparison, the maximal excess method does this
for $\,n=5$, $9$, and $17$, but not for $\,n=13$ or $21$. For
$\,n=45$, $53$, and $69$, our method reproduces the largest
known determinants, reported in \cite{OSDS}. Most significantly,
we have used it to set new records for $\,n=49$, $65$, $73$,
$77$, $85$, $93$, and $97$.
%
%

Theorem \ref{thm:construction} makes it clear that to set records
of this type using our construction, we need to find
{3-normalized} Hadamard matrices of large excess. In
particular, to improve upon the maximal excess
technique, we need to satisfy condition (\ref{eqn:compare}).

The general literature of Hadamard matrices provides a variety of
techniques for constructing Hadamard matrices, and explicit
examples of relatively low order ($\,n<200\,$, say) abound on
the internet. Unfortunately, we are unable to predict, based on
construction technique or structural elements, for example,
which of these matrices will have large excess after
3-normalization.

We have therefore searched for large-excess 3-normalized Hadamard
matrices using an approach with two main steps. Namely,%

\begin{enumerate}

\item
Generate initial Hadamard matrices from distinct equivalence
classes in $\,\had{n}\,$, and then

\item
Find a matrix in each class which attains the maximum excess
subject to the 3-normalization constraint.
\end{enumerate}

A program we wrote to perform the latter task runs in polynomial
time and has served us well for $\,n<100\,$ or so. Finding
initial Hadamards which can be 3-normalized with high excess,
however, poses real difficulty in general. Starting in dimension
$\,56\,$ for the case $\,n\equiv0\,$, and dimension
92 for the case $\,n\equiv 4\pmod8\,$, we have not yet succeeded
in finding 3-normalized Hadamard matrices with excess achieving
our upper bound $\,\nu_n^*\,$.

We next discuss our approach to the two tasks above.

\subsection*{Generating initial Hadamards}
For orders $n\equiv 4\pmod{8}$, online libraries of Williamson
and so-called ``good'' matrices \cite{Se} provide a rich source
of initial Had\-amard matrices. When $n\equiv 0 \pmod{8}$,
however, one finds fewer examples online, and we have resorted
to generating our own from smaller Hadamards. The tensor product
of Hadamard matrices $\,H_1\in\had{s_1}\,$ and
$\,H_2\in\had{s_2}\,$ always belongs to $\,\had{s_1s_2}\,$. More
generally, one can construct Hadamard matrices of order
$\,s_1s_2\,$ by a technique known as ``weaving,'' given
$\,s_1\,$ matrices in $\,\had{s_2}\,$ and $\,s_2\,$ matrices in
$\,\had{s_1}\,$ \cite{Cr}. Craigen and Kharaghani have ``woven''
Hadamard matrices with very high excess \cite{CK}. Unfortunately
weaving does not, in general, produce \textit{3-normalized}
Hadamard matrices of high excess.

On the other hand, we have used another technique, the
Multiplication Theorem of Agaian-Sarukhanyan, with quite good
results. This Theorem constructs a Hadamard of order
$\,{8k_1k_2}\,$, given $\,H_1\in\had{4k_1}\,$ and
$\,H_2\in\had{4k_2}\,$. The reader may consult \cite{Ag} or
\cite{SY} for further details.

\subsection*{Excess-maximization}\label{ssec:algorithm}%
Our algorithm for maximizing excess on an equivalence class of
3-normalized Hadamard matrices relies on the fact that the
excess depends only on a choice of the three rows used in the
3-normalization. We make this fact precise in Lemma
\ref{lem:algorithm} below, which requires the following
definition:

\begin{definition}\label{def:spm}
Let $\,S_n^\pm\,$ denote the group of $\,\nxn\,$ \textit{signed
permutation matrices}, i.e.~matrices obtained by permuting the
rows of an $\,\nxn\,$ diagonal matrix of $\,\pm1$s. Let
$\,S_{3,n-3}^\pm\subset S_n^\pm\,$ denote the subgroup
comprising matrices with $\,(3,n-3)\,$ block-diagonal form.
\end{definition}

\begin{lem}\label{lem:algorithm}
Suppose $\,N\,$ is a 3-normalized Hadamard matrix of order
$\,n>2\,$. Then any 3-normalized Hadamard matrix $\,N'\,$
equivalent to $\,N\,$ can be written as
\begin{equation}\label{eqn:equiv}
N' =R\,N\, C\ ,\qquad\text{where $\,R,\,C\in S_n^\pm\,$}\ ,
\end{equation}
and $\,\ex{N'}\,$ depends only on the coset to which
$\,R\,$ belongs in $\,S_n^\pm/S_{3,n-3}^\pm\,$.
\end{lem}

\begin{proof}
Equation (\ref{eqn:equiv}) simply restates the definition of
``equivalent'' given in Defn.~\ref{def:equiv}. The main
conclusion here is the one about excess, and for that, it
suffices to show $\,\ex{RNC}=\ex{N}\,$ whenever both matrices
are 3-normalized, and $\,R\in S_{3,n-3}^\pm\,$.

To do so, observe that rows 1,2, and 3 of $\,N\,$
form a $3\times n$ submatrix with columns
\begin{equation*}
\left[\begin{array}{c} +\\+\\+\end{array}\right]\,,\quad
\left[\begin{array}{c} -\\-\\+\end{array}\right]\,,\quad
\left[\begin{array}{c} -\\+\\-\end{array}\right]\,,\ \text{and}\
\left[\begin{array}{c} +\\-\\-\end{array}\right]\ ,
\end{equation*}
where $\,+,\,-\,$ signify $\,\pm1\,$ respectively. These
vectors occur in successive blocks of length $\,k:=n/4\,$, no
two are mutually opposite, and signed row permutations preserve
their distinctness. It follows that the signed permutation
$\,C\,$ in (\ref{eqn:equiv}) can permute columns within any
one of the four $k$-blocks, or move entire $k$-blocks,
but it must maintain the contiguity of each block. Further,
one sees that $\,C\,$ negates all columns or none, depending on
whether $\,R\,$ negates an odd or even number, respectively, of
the first three rows.

Now recall from Proposition \ref{prop:nor} that for each
$\,i=4,5,\dots,n\,$, row $\,i\,$ of $\,N\,$ has the form
$\,(\va_i\,,\,\vb_i\,,\,\vc_i\,,\,\vd_i)\,$, where
$\,\va_i,\,\vb_i,\,\vc_i,\,\vd_i\in\{\pm1\}^k\,$ all have equal
excess. The facts above now force row $\,i\,$ of $\,NC\,$ to
take the form
\[
(\va'_i\,,\,\vb'_i\,,\,\vc'_i\,,\,\vd'_i)\ ,
\]
where we get $\,\{\va'_i,\,\vb'_i,\,\vc'_i,\,\vd'_i\}\,$ by
permuting $\,\{\va_i,\,\vb_i,\,\vc_i,\,\vd_i\}\,$, permuting the
entries of each vector, and then, possibly, negating them all.
Moreover, since $\,R\,$ belongs to the subgroup
$\,S^\pm_{3,n-3}\,$ the last $\,n-3\,$ rows of $\,R\cdot NC\,$
equal those $\,NC\,$, modulo permutations and sign changes of
entire rows. It follows that for some permutation $\,\sigma\,$ of
$\,\{4,5,\dots,n\}\,$, we have
\begin{equation}\label{eqn:pmex}
\ex{\va'_i} = \ex{\vb'_i} = \ex{\vc'_i} = \ex{\vd'_i}
= \pm\ex{\va_{\sigma(i)}}\ ,
\end{equation}

But the excesses of rows 4 through $n$ of $\,RNC\,$ must be
non-negative (Definition \ref{def:nor}), so the final sign in
(\ref{eqn:pmex}) must always be ``$+$''. That is, rows
4 through $n$ of $\,N\,$ and $\,RNC\,$ have exactly the same
excesses, modulo permutation. These rows account for the entire
excess of a 3-normalized Hadamard matrix, so we are done.
\end{proof}

The subgroup $\,S^\pm_{3,n-3}\subset S^\pm_3\,$ is clearly isomorphic to
$\,S_3^\pm\times S_{n-3}^\pm\,$, and the order of $\,S_n^\pm\,$
is clearly $\,2^nn!\,$ for any $\,n\,$.  We therefore have
\begin{equation}\label{eqn:qs}
\#\left(S^\pm_n/S_{3,n-3}^\pm\right) =
{n\choose 3}\ .
\end{equation}
Indeed, the $\,n\choose 3\,$ matrices we get by choosing three
rows of the $\nxn$ identity matrix, then exchanging them
(in any order) with rows 1, 2, and 3, each represent a
different coset of $\,S^\pm_n/S^\pm_{3,n-3}\,$.

Using these facts we easily find a 3-normalized matrix of
maximal excess equivalent to any initial $\,H\in\had{n}\,$.
We take $\,\mu:={n\choose3}\,$
matrices $\,\left\{R_\a\right\}_{\a=1}^\mu\subset S^\pm_n\,$,
each representing a different coset of
$\,S^\pm_n/S^\pm_{3,n-3}\,$, and then form the matrices
\begin{equation}\label{eqn:rh}
R_1H,\,R_2H,\,\dots,\,R_\mu H\in\had{n}
\end{equation}
According to Proposition \ref{prop:nor}, we can 3-normalize each
$\,R_iH\,$ using only a signed column permutation $\,C_i\,$, and
some row-negations, which we can represent by a diagonal matrix
$\,D_i\,$ of $\,\pm1$s, multiplying $\,R_i\,$ on the right. Such
diagonal matrices clearly belong to $\,S^\pm_{3,n-3}\,$, so
we now have $\,{n\choose3}\,$ 3-normalized Hadamard matrices
$\,N_i\,$, $\,i=1,2,\dots,{n\choose 3}\,$, with
\[
N_i = R_iD_iHC_i\ .
\]
Since the $\,R_iD_i$s represent every coset of
$\,S^\pm_n/S^\pm_{3,n-3}\,$, Lemma \ref{lem:algorithm} above
guarantees that among these, the $\,N_i\,$ with largest excess
actually maximizes excess among all 3-normalized Hadamard matrices
equivalent to $\,H\,$.

\begin{remark}
Already when $\,n=16\,$, one of the five equivalence
classes in $\,\had{16}\,$ fails to contain any 3-normalized matrix
that achieves the bound in Theorem \ref{thm:bounds}.

By tabulating the row sums one gets after 3-normalizing
representatives from each of the $\,{n\choose 3}\,$ cosets
described above, one associates to any $\,H\in\had{n}\,$ a
statistic equivalent to the \textit{profile} defined by Cooper,
Milas, \& Wallis in \cite{CMW}.
\end{remark}

\subsection*{A table of large determinants}\label{ssec:table}
For each $\,n\equiv0\pmod4\,$, $\,4\le n\le 100\,$, Table 1 below
displays the results of our efforts to find 3-normalized
Hadamard matrices of large excess, and the large determinants
they generate via our construction. Additionally, the table
compares our results to those obtained using the maximal excess
technique of Farmakis \& Kounias, and indicates the largest
known determinant of each order. The actual matrices which
attain these determinants are posted on our website \cite{OS}.

We now give a detailed explanation of Table 1.

\begin{table}\label{table:main}
\begin{tabular}{c|rr|rr|rrc}
  & \multicolumn{2}{c}{\emph{Upper bounds}} & \multicolumn{2}{c}{\emph{Best known}} & \multicolumn{3}{c}{}\\
$n+1$ &  $n+\nu_{n}^*$ & $\sigma_{n}^*$ & $n+\nu_{n}$ & $\sigma_{n}$ & $\mu_{n+1}$ & $\beta_{n+1}$ & $\mu_{n+1}/\beta_{n+1}$\\ \hline\hline
      &            &       &         &         &\ph{$\ \ddag1234$}& \ph{$12345$}\\
  5   & \bf 8     &\bf  8 &\bf 8    &\bf  8   & $\underline{3}^*\phantom{\sf5/7}$      & 3\phantom{.00}  & 1\ph{.00}  \\
  9   & 16        &\bf 20 & 16      &\bf 20   & $7^*\phantom{\sf5/7}$      & 8.25   & 0.85 \\
  13  & {\bf 48}  & 36    & \bf 48  & 36      & $\underline{15}^*\phantom{\sf5/7}$     & 15\phantom{.00}  & 1\ph{.00} \\
  17  & 48        &\bf 64 & 48      &\bf 64   & $20^*\phantom{\sf5/7}$     & 22.98  & 0.87 \\
  21  & {\bf 96}  & 80    & \bf 96  & 80      & $\underline{29}^*\phantom{\sf5/7}$     & 32.02  & 0.91 \\
  25  & 96        &\bf 112& 96      & 112     & $42^*\phantom{\sf5/7}$  & 42\phantom{.00}  & 1\ph{.00} \\
  29  & {\bf 152} & 140   & 152     & 140     & $45\sf5/7\phantom{{}^*}$ & 52.85 & 0.87 \\
  33  & 160       &\bf 172& 160     & 172     & $55\sf1/8\phantom{{}^*}$ & 64.50 & 0.85 \\
  37  & \bf 216   &\bf 216& 216     & 216     & ${\bf72}\phantom{\sf5/7^*}$   & 76.90  & 0.94 \\
  41  &  240      &\bf 244& 240     & 244     & $90^*\phantom{\sf5/7}$  & 90\phantom{.00}  & 1\ph{.00} \\
  45  & {\bf 288} & 280   &\bf 288  & 280     & $\underline{83}\phantom{\sf5/7^*}$  & 103.77 & 0.80 \\
  49  & {\bf 336} & 324   &\bf 336  & 324     & $\underline{\bf96}\phantom{\sf5/7^*}$  & 118.19 & 0.81 \\
  53  & {\bf 368} & 364   &\bf 368  & 364     & $\underline{105}\phantom{\sf5/7^*}$ & 133.21 & 0.79 \\
  57  & {\bf 448} & 408   & 384     & \bf 400 & $114\phantom{\sf5/7^*}$ & 148.82 & 0.77 \\
  61  & {\bf 456} & 452   & 456     & 440     & $165^*\phantom{\sf5/7}$ & 165\phantom{.00} & 1\ph{.00} \\
  65  & {\bf 560} & 512   &\bf 528  & 512     & $\underline{\bf 148}\phantom{\sf5/7^*}$ & 181.73 & 0.81 \\
  69  & {\bf 552} & 548   &\bf 552  & 544     & $\underline{155}\phantom{\sf5/7^*}$ & 198.98 & 0.78 \\
  73  & {\bf 656} & 600   &\bf 624  & 580     & $\underline{\bf174}\phantom{\sf5/7^*}$ & 216.75 & 0.80 \\
  77  & {\bf 656} & 652   &\bf 656  & 628     & $\underline{\bf183}\phantom{\sf5/7^*}$ & 235.02 & 0.78 \\
  81  & {\bf 752} & 704   &\bf 704  &\bf 704  & $\underline{196}\phantom{\sf5/7^*}$ & 253.77 & 0.78 \\
  85  & {\bf 768} & 756   &\bf 768  & 756     & $\underline{\bf213}\phantom{\sf5/7^*}$ & 273\phantom{.00} & 0.78 \\
  89  & {\bf 864} & 812   &768     &\bf 792  & $220\phantom{\sf5/7^*}$ & 292.69 & 0.75 \\
  93  & {\bf 888} & 872   &\bf 864  & 828     & $\underline{\bf239}\phantom{\sf5/7^*}$ & 312.83 & 0.76 \\
  97  & {\bf 976} & 932   &\bf 928  & 920     & $\underline{\bf256}\phantom{\sf5/7^*}$ & 333.42 & 0.77 \\
 101  & {\bf 1016}& 1000  & 984     &\bf 1000 & $275\phantom{\sf5/7^*}$ & 354.44 & 0.78 \\

\end{tabular}
\caption{\strut Largest known excesses of order
$\,n\equiv0\pmod4\,$ and determinants of order $\,n+1\,$.
$\,\mu_{n+1}\,$ encodes the largest known determinant of order
$\,n+1\,$. Compare $\,\mu_{n+1}\,$ with the corresponding factor
for our determinant ($\,{(n+\nu_n)+n\over4}\,$), and that of
\cite{FK} ($\,{\sigma_n+n\over4}\,$).}
\end{table}

\subsubsection*{Notation}

\begin{itemize}

\item
$\,\s_{n}^*\,$ and $\,\s_n\,$ respectively denote the sharpest
known upper bound for the maximal excess, and the largest known
excess, of a Hadamard matrix of order $\,n\,$.

\item
$\,\nu_n^*\,$ and $\,\nu_n\,$ respectively denote the sharpest
known upper bound for the excess of, and the largest known excess
for, a 3-normalized Hadamard matrix of order $\,n\,$.

\item
$\,\mu_{n+1}\,$ and $\,\beta_{n+1}\,$ abbreviate the largest
known determinant $\det M_{n+1}$, and the Barba bound
$\,B(n+1)\,$ (see (\ref{eqn:barba})), by dividing out large
integer factors as follows:
\begin{equation*}\label{eqn:mu&beta}
\mu_{n+1}:= {\det M_{n+1}\over 2^n\cdot k^{2k-1}}\ ,\qquad
\beta_{n+1}:= {B(n+1)\over 2^n\cdot k^{2k-1}}\ .
\end{equation*}
Here $\,k=n/4\,$, and $\,M_{n+1}\,$ is a sign matrix of order
$\,n+1\,$ having largest known determinant.

\end{itemize}
\goodbreak

\subsubsection*{Explanatory notes}

\begin{enumerate}

\item\textit{The ``Upper bounds'' columns.}
The tabulated values for $\,\s_n^*\,$ all come from \cite[Table
2]{FK}.  Except for $\,n=80\,$ (see Remark \ref{rem:n=80}), we
computed the values for $\,\nu_n^*\,$ using Theorem
\ref{thm:bounds}, and we list $\,n+\nu_n^*\,$ instead of just
$\,\nu_n^*\,$, to facilitate the key comparison indicated by
(\ref{eqn:compare}). In each row, boldface indicates the
larger of $\,n+\nu_n^*\,$ and $\,\s_n^*\,$.
\medskip

\item\textit{The ``Best known'' columns.}
A bold entry in either column indicates that the largest known
determinant of that order is produced by applying the
corresponding construction to a matrix with the associated
excess $\,\nu_n\,$ or $\,\sigma_n\,$. When neither entry
is bold, neither method constructs the current record.

\noindent
All values for $\,\s_{n}\,$ come from  \cite[Table 2]{FK}
except for $\,\sigma_{72}=580\,$ and $\,\sigma_{76}=628\,$,
which were reported in \cite{OSDS}.
\medskip

\item\textit{The ``$\mu_{n+1}$'' column.}
Current records which can be constructed by our method are
underlined. A bold entry represents a determinant
that achieves a new record, reported for the first time here.
Values of $\,\mu_{n+1}\,$ corresponding to determinants known to
attain $\,\md(n+1)\,$ are marked with an asterisk.
\medskip

\item
For $\,n\le 20\,$, either the 3-normalized or maximal excess
method (or both) constructs an order $\,n+1\,$ sign matrix that
achieves the global determinant maximum.
Equivalent matrices were all constructed by earlier
investigators, however: See \cite{Mo,Wi} for $\,n+1=5\,$,
\cite{EZ} for $\,n+1=9\,$, \cite{R} for $\,n+1=13\,$, \cite{MK}
for $\,n+1=17\,$, and \cite{CKM} for $\,n+1=21\,$.
\medskip

\item
For sizes $\,n+1=25, 29, 33, 37, 41\,$ and $\,61\,$, neither our
method nor that of Farmakis \& Kounias constructs the largest
known determinant. In dimensions $\,25, 41$, and
$61$, matrices coming from block designs (SBIBDs) attain 100\% of
Barba's bound $\,B(n+1)\,$; see \cite{Bh, R} for $\,n+1=25\,$,
\cite{BHH} for $\,n+1=41\,$, and \cite{Br} for $\,n+1=61\,$. The
largest known determinant for $\,n+1=33\,$ was reported in
\cite{OSDS}. For $\,n+1=29$ and $37\,$, we sketch a construction
for the largest known determinants below.

\end{enumerate}

\section{Improved results for $n+1=29,\,37$}\label{sec:symmetrize}

The ``Upper bounds'' and ``Best known'' columns of Table 1 show
that for $n+1=29$, our construction achieves a determinant
larger than any the maximal excess method can produce. The value
$\,\mu_{29}=45\sf5/7\,$ there, however, corresponds to the still
larger determinant reported in~\cite{OSDS} using a gradient
ascent algorithm. We have since discovered that the latter
matrix can be formed by applying a symmetric version of our
construction. To do this, we start with a matrix
$\,H\in\had{28}\,$ for which both $\,H\,$ and $\,H^T\,$ satisfy
a generalized 3-normalization condition, and for which both have
the optimal set of row sums prescribed by
Theorem~\ref{thm:bounds}. We then apply appropriately adapted
versions of steps 1--4 of our construction (\S
\ref{sec:construction}) to $\,H\,$, transpose the result,
and apply the same four steps again. This produces a matrix
equivalent to the one reported in \cite{OSDS}.

Hoping to profit by applying the same idea to other orders, we
have so far achieved one notable success: A matrix
$\,A\in\sm{37}\,$ with $\det A = 72\times9^{17}\times2^{36}$.
This determinant, reported here for the first time, attains 94\%
of the Barba bound, and we conjecture it to be the global maximum
for its order. By way of comparison, Table 1 
shows that the largest determinant that either the maximal excess
method or our own can hope produce is
$63\times9^{17}\times2^{36}$. This value has been
achieved by both methods, but it represents just 82\% of Barba's
bound.

Guided by the $\,n+1=29\,$ case, we start by finding a matrix
$\,H\in\had{36}\,$ which can be ``doubly 3-normalized'' with
optimal row and column sums. In this dimension, the three
3-normalizing rows and columns overlap in a somewhat more
complicated way than in dimension 28. This forces
further modifications of the double 4-step construction
described above, but the procedure is very similar.

We are currently investigating the relationship between these two
examples and possible generalizations. We will publish further
details of this work in our sequel to the present work.
Meanwhile, we have posted both the $\,n+1=29\,$ and $\,n+1=37\,$
examples on our website \cite{OS}.

\section*{Acknowledgements}
Our work was greatly facilitated by Hadamard matrices provided by
the websites of N.J.A. Sloane, J. Seberry, C. Koukouvinos, and E.
Spence. We gratefully acknowledge them for making these
available. We used \emph{Mathematica} extensively during this
project, and we also thank Indiana University for its suppport,
and for the use of its Sun E10000 computing platform.
\bigskip


\bigskip

\end{document}